# Tareas matemáticas abiertas como respuesta didáctica a la inteligencia artificial generativa en contextos post-IA

# Open Mathematical Tasks as a Didactic Response to Generative Artificial Intelligence in Post-AI Contexts

**Felix De la Cruz Serrano**.
I.E. María Reiche.
Perú.
feldese@gmail.com

## Resumen

La disponibilidad generalizada de herramientas de inteligencia artificial generativa plantea nuevos desafíos para la enseñanza de las matemáticas escolares, particularmente en relación con el sentido formativo de las tareas tradicionales. En contextos post-IA, numerosas actividades pueden resolverse de manera automática, sin que los estudiantes se involucren en procesos de interpretación, toma de decisiones o validación matemática.

Este estudio analiza una experiencia de aula en educación secundaria en la que se implementan tareas matemáticas abiertas como respuesta didáctica a este escenario, con el propósito de sostener la actividad matemática del estudiante. Desde un enfoque cualitativo de carácter descriptivo–interpretativo, se examinan las formas de trabajo matemático que emergen durante la resolución de las tareas, mediadas por el dispositivo didáctico de regulación COMPAS.

El análisis se organiza en torno a cuatro ejes: el diseño de tareas abiertas en contextos post-IA, la agencia matemática del estudiante, la complementariedad humano–IA y las prácticas de modelización y validación. Los resultados muestran que, bajo una regulación didáctica explícita, los estudiantes mantienen el control epistémico de la actividad matemática, incluso en presencia de inteligencia artificial generativa.

*Palabras clave:* tareas matemáticas abiertas; inteligencia artificial generativa; educación matemática; complementariedad humano–IA; control epistémico.

## Abstract

The widespread availability of generative artificial intelligence tools poses new challenges for school mathematics education, particularly regarding the formative role of traditional mathematical tasks. In post-AI educational contexts, many activities can be solved automatically, without engaging students in interpretation, decision-making, or mathematical validation processes.

This study analyzes a secondary school classroom experience in which open mathematical tasks are implemented as a didactic response to this scenario, aiming to sustain students' mathematical activity. Adopting a qualitative and descriptive–interpretative approach, the study examines the forms of mathematical work that emerge during task resolution, mediated by the didactic regulation device COMPAS.

The analysis is structured around four analytical axes: open task design in post-AI contexts, students' mathematical agency, human–AI complementarity, and modeling and validation practices. The findings suggest that, under explicit didactic regulation, students retain epistemic control over mathematical activity, even in the presence of generative artificial intelligence.

*Keywords:* open mathematical tasks; generative artificial intelligence; mathematics education; human–AI complementarity; epistemic control.

**Introducción**

La creciente disponibilidad de herramientas de inteligencia artificial generativa está transformando de manera significativa las condiciones en las que se desarrolla la enseñanza y el aprendizaje de las matemáticas escolares. Estas herramientas, capaces de producir cálculos, procedimientos y explicaciones matemáticas de forma inmediata, se encuentran hoy al alcance de estudiantes y docentes tanto dentro como fuera del aula, lo que plantea nuevos desafíos para las prácticas educativas tradicionales (Hemmer et al., 2024).

En particular, la presencia de la inteligencia artificial tensiona el sentido formativo de muchas tareas matemáticas escolares, ya que numerosas actividades pueden resolverse sin que los estudiantes se involucren en procesos de interpretación, toma de decisiones o validación. Desde la didáctica de la matemática, se ha advertido que este desplazamiento refuerza una visión pragmática de la actividad matemática, centrada en la obtención de resultados, en detrimento de su dimensión epistémica y de las prácticas de pensamiento que le otorgan significado (Artigue, 2025; Schoenfeld, 2020).

En este trabajo, el término post-IA se utiliza para referirse a contextos educativos en los que la inteligencia artificial generativa se encuentra disponible de manera generalizada. En dichos contextos, el diseño de tareas y la regulación de la actividad matemática deben asumir la presencia de la IA como una condición del entorno, y no como un recurso externo, excepcional o susceptible de ser simplemente prohibido.

Frente a este escenario, las respuestas centradas exclusivamente en la restricción del uso de la tecnología o en su incorporación acrítica como herramienta de apoyo resultan insuficientes. El problema central no reside en la herramienta en sí, sino en el tipo de actividad matemática que las tareas escolares promueven y en las decisiones didácticas que regulan su uso. Esta actividad se entiende como el conjunto de procesos de formulación, razonamiento y validación que caracterizan el trabajo matemático escolar (Niss & Højgaard, 2019). En este sentido, el diseño de tareas emerge como un elemento clave para sostener la participación activa del estudiante en contextos de disponibilidad generalizada de tecnologías avanzadas.

Dentro de este marco, las tareas matemáticas abiertas —caracterizadas por la indeterminación inicial, la necesidad de formular supuestos y la exigencia de validar resultados— han sido identificadas como una respuesta didáctica pertinente para promover la agencia del estudiante y el control epistémico sobre la actividad matemática (Sullivan et al., 2013; Leung & Baccaglini-Frank, 2017). No obstante, su potencial en contextos con inteligencia artificial depende de la existencia de dispositivos de regulación que orienten la interacción entre estudiantes, tareas y tecnología.

Este estudio analiza una experiencia de aula en educación secundaria en la que se implementan tareas matemáticas abiertas mediadas por un dispositivo didáctico de regulación, con el propósito de examinar las formas de trabajo matemático que emergen en contextos de disponibilidad de inteligencia artificial generativa. Desde un enfoque cualitativo de carácter descriptivo–interpretativo, se analizan las prácticas desarrolladas por los estudiantes atendiendo al diseño de las tareas, a su agencia matemática y a los procesos de modelización y validación implicados.

Los resultados sugieren que, cuando el uso de la inteligencia artificial es regulado mediante decisiones didácticas explícitas, los estudiantes mantienen el control epistémico de la actividad matemática, asumiendo un rol activo en la construcción y validación de modelos. Estos hallazgos aportan elementos relevantes para el debate actual sobre la enseñanza de las matemáticas en contextos post-IA y subrayan el papel del diseño de tareas y de la regulación didáctica en la preservación del sentido formativo de la actividad matemática.

## Marco teórico

### Complementariedad humano–IA en matemáticas escolares

El análisis del estudio se apoya en la noción de **complementariedad humano–IA**, entendida como una forma de colaboración en la que las capacidades humanas y las de los sistemas de inteligencia artificial no se sustituyen, sino que se distribuyen de manera asimétrica y dinámica según la naturaleza de la tarea (Hemmer et al., 2024). Desde esta perspectiva, la inteligencia artificial puede asumir funciones asociadas a la exploración rápida de resultados, la generación de ejemplos o el apoyo instrumental, mientras que los estudiantes conservan el control sobre la interpretación de la situación, la toma de decisiones y la validación matemática.

En el contexto de las matemáticas escolares, esta complementariedad adquiere un carácter eminentemente didáctico. La presencia de IA generativa introduce el riesgo de que la actividad matemática del estudiante se reduzca a la aceptación de respuestas producidas externamente. Sin embargo, cuando el diseño de la tarea exige formular supuestos, justificar elecciones y evaluar la coherencia de los resultados, la IA se integra como un recurso subordinado a la actividad intelectual del estudiante, preservando el control epistémico del proceso matemático (Hemmer et al., 2024).

### Diseño de tareas abiertas en contextos con inteligencia artificial

El estudio se sitúa en el marco del diseño de **tareas matemáticas abiertas**, caracterizadas por presentar situaciones con indeterminación inicial, ausencia de un único procedimiento de resolución y necesidad de validar los resultados obtenidos. Investigaciones previas han mostrado que este tipo de tareas favorece la agencia del estudiante y promueve prácticas matemáticas centradas en la exploración, la argumentación y la toma de decisiones (Sullivan et al., 2013; Leung & Baccaglini-Frank, 2017).

En contextos de disponibilidad de inteligencia artificial generativa, las tareas abiertas adquieren un papel particularmente relevante. Frente a tareas cerradas y algorítmicas, fácilmente resolubles mediante IA, las tareas abiertas dificultan una delegación automática de la actividad matemática y obligan a los estudiantes a involucrarse en procesos de interpretación y modelización. De este modo, el diseño de la tarea se convierte en un elemento central para sostener la dimensión epistémica del trabajo matemático y para regular el uso de la tecnología en el aula.

### COMPAS como dispositivo didáctico de regulación

Para articular la complementariedad humano–IA en la práctica de aula, el estudio introduce el dispositivo didáctico **COMPAS**, concebido como una estructura de regulación de la actividad matemática en contextos con inteligencia artificial generativa. COMPAS no se plantea como un método prescriptivo, sino como un **dispositivo conceptual** que orienta la gestión de la tarea y la interacción con la tecnología a través de seis fases articuladas: Comprensión de la situación, Organización de supuestos y datos, Modelización matemática, Producción y exploración de resultados, Análisis y validación, y Socialización y síntesis.

Estas fases no constituyen una secuencia rígida, sino un marco flexible que permite regular los momentos de uso de la IA y favorecer la participación activa de los estudiantes. En particular, COMPAS busca evitar la clausura prematura de la tarea y promover prácticas de validación colectiva, en coherencia con enfoques didácticos que destacan la importancia de la discusión y la argumentación en la construcción del conocimiento matemático (Stylianides & Stylianides, 2019).

## Metodología

### Enfoque metodológico y ejes de análisis

El estudio adopta un enfoque cualitativo, de carácter descriptivo–interpretativo, orientado al análisis de una experiencia de aula desarrollada en educación secundaria. El interés se centra en comprender las formas de trabajo matemático que emergen durante la resolución de tareas abiertas en contextos de disponibilidad de inteligencia artificial generativa, más que en medir efectos o comparar desempeños.

La experiencia se implementó en sesiones regulares de clase, integradas al desarrollo habitual del curso de matemáticas. La recolección de datos incluyó producciones escritas de los estudiantes, registros de interacciones en el aula y observaciones del docente-investigador.

El análisis se organizó en torno a cuatro ejes analíticos, derivados del marco teórico y de las decisiones didácticas adoptadas:

(a) El diseño de tareas abiertas en contextos con IA.
(b) La agencia matemática del estudiante.
(c) La complementariedad humano–IA.
(d) Las prácticas de modelización y validación.

Estos ejes estructuran la presentación y el análisis de los resultados del estudio.

## Desarrollo y resultados

En esta sección se presenta el análisis de las formas de trabajo matemático que emergen durante la resolución de tareas abiertas en contextos de disponibilidad de inteligencia artificial generativa. El análisis se organiza en torno a cuatro ejes, definidos a partir del marco teórico y metodológico del estudio, que permiten dar cuenta de la actividad matemática desarrollada y de las tensiones asociadas al uso de la IA.

### Diseño de tareas abiertas post-IA y configuración de la actividad matemática

El primer eje de análisis se centra en el papel del diseño de las tareas abiertas en la configuración de la actividad matemática desarrollada por los estudiantes en contextos de disponibilidad de inteligencia artificial generativa. En particular, se analizan aquellas decisiones de diseño que dificultan una resolución automática y obligan a asumir decisiones matemáticas significativas desde las etapas iniciales de la tarea.

Las tareas implementadas se caracterizan por la ausencia deliberada de datos numéricos explícitos y de consignas cerradas, lo que genera un espacio inicial de indeterminación. Esta indeterminación no constituye una falta de información, sino una decisión de diseño orientada a desplazar el foco de la actividad desde el cálculo hacia la interpretación de la situación y la formulación de supuestos relevantes. Un ejemplo representativo de este tipo de diseño es la tarea ¿Cuántas botellas se necesitan?, presentada a partir de una situación cotidiana sin información cuantitativa explícita (Figura 1).

La imposibilidad de realizar un cálculo directo impide una delegación inmediata de la resolución en herramientas de inteligencia artificial y sitúa la actividad matemática en un plano previo a cualquier procesamiento algorítmico. En este contexto, los estudiantes deben decidir qué variables considerar, qué supuestos adoptar y qué relaciones establecer entre ellas, configurando así el núcleo de la actividad matemática.

El análisis de las producciones y de las interacciones muestra que este tipo de diseño promueve una actividad matemática orientada a la construcción de sentido, en la que los resultados numéricos adquieren significado únicamente en función de los supuestos que los sustentan. Desde esta perspectiva, el diseño de la tarea actúa como un mecanismo de regulación epistémica, al impedir la clausura temprana de la actividad y al subordinar cualquier uso de la inteligencia artificial a decisiones matemáticas previas.

En contextos post-IA, este tipo de diseño resulta particularmente relevante, ya que redefine la función de la inteligencia artificial dentro del proceso de resolución. La IA puede ser utilizada como un recurso auxiliar para explorar o contrastar resultados, pero no como una fuente de soluciones definitivas. De este modo, el diseño de tareas abiertas contribuye a preservar la dimensión epistémica de la actividad matemática y a evitar que la resolución se reduzca a la aceptación acrítica de respuestas producidas externamente.

**Figura 1**

*Situación presentada a los estudiantes para la tarea «¿Cuántas botellas se necesitan?». La imagen se introdujo sin datos numéricos explícitos ni consignas cerradas, con el propósito de favorecer la interpretación de la situación y la formulación de supuestos.*

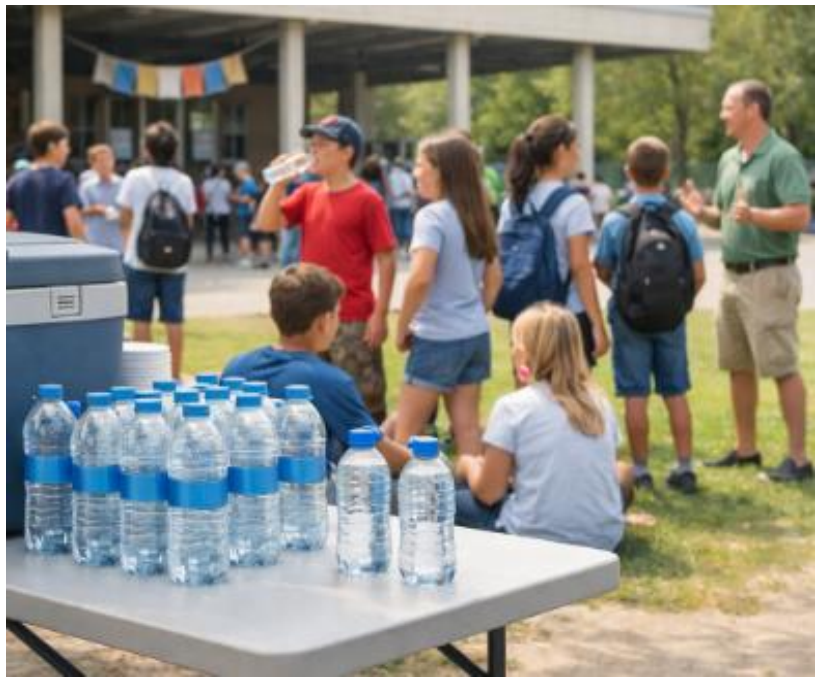

**Agencia matemática del estudiante en la resolución de tareas con IA**

El segundo eje de análisis se centra en la agencia matemática del estudiante, entendida como la capacidad de asumir control sobre las decisiones que estructuran el proceso de resolución, más allá de la mera producción de respuestas. En el contexto de tareas abiertas mediadas por inteligencia artificial generativa, esta agencia se manifiesta principalmente en la formulación, revisión y negociación de los supuestos que sustentan los modelos matemáticos construidos.

El análisis de las producciones estudiantiles y de las interacciones desarrolladas durante la experiencia muestra que los estudiantes asumen un rol activo en la resolución de las tareas, particularmente cuando deben decidir qué variables considerar, qué valores asignar y qué relaciones establecer entre ellas. La necesidad de determinar, por ejemplo, el número de participantes o el consumo esperado por persona, favorece la emergencia de diferentes estrategias y posicionamientos frente a la tarea, que no pueden ser resueltos de manera automática ni externalizados completamente a la tecnología.

Desde esta perspectiva, la agencia matemática no se expresa únicamente en la obtención de resultados, sino en la capacidad de tomar decisiones, justificar elecciones y cuestionar la pertinencia de los supuestos adoptados. Estas prácticas se hacen visibles en las discusiones colectivas, en las que los estudiantes comparan modelos, evalúan su coherencia con la situación planteada y revisan críticamente sus planteamientos iniciales. La actividad matemática se configura, así como un proceso de construcción progresiva, en el que los modelos son susceptibles de ser ajustados y reformulados.

Asimismo, cuando los estudiantes recurren a la inteligencia artificial como recurso de apoyo, el análisis muestra que las decisiones matemáticas centrales permanecen bajo su control. Las respuestas o cálculos generados por la IA son sometidos a contraste y validación en función de los supuestos previamente establecidos, en lugar de ser aceptados de manera acrítica. De este modo, la inteligencia artificial no sustituye la actividad matemática del estudiante, sino que se integra como un recurso subordinado dentro de un proceso de resolución regulado por decisiones humanas.

En este contexto, la agencia matemática del estudiante se ve reforzada por la necesidad de evaluar la plausibilidad y coherencia de las soluciones propuestas, incluso cuando estas provienen de una fuente tecnológicamente avanzada. La actividad matemática se configura así como un espacio de reflexión, argumentación y validación, en el que el estudiante mantiene un rol activo y responsable, preservando el control epistémico del proceso de resolución en contextos de disponibilidad de inteligencia artificial generativa.

**Complementariedad humano–IA en la práctica de aula**

El tercer eje de análisis aborda la complementariedad humano–IA como una forma específica de distribución de funciones durante la resolución de tareas matemáticas abiertas en contextos de disponibilidad de inteligencia artificial generativa. El análisis de la experiencia pone de manifiesto una distribución asimétrica de roles, en la que las capacidades de la tecnología y las de los estudiantes no se superponen, sino que se articulan de manera funcional en función de las exigencias de la tarea.

En la práctica observada, la inteligencia artificial es utilizada principalmente como un recurso de apoyo instrumental, orientado a la exploración de resultados posibles, la realización de cálculos puntuales o el contraste de valores obtenidos. En contraste, los estudiantes conservan la responsabilidad de interpretar la situación, formular supuestos y decidir la pertinencia de las soluciones propuestas. Esta distribución de funciones sitúa las decisiones matemáticas centrales fuera del alcance de una automatización completa.

Esta forma de complementariedad no se produce de manera espontánea, sino que es el resultado de un diseño de tareas que exige decisiones matemáticas previas a cualquier cálculo y de una regulación didáctica explícita que evita la clausura temprana de la actividad. En este contexto, la inteligencia artificial no actúa como una fuente de respuestas definitivas, sino como un recurso parcial cuya información debe ser evaluada, interpretada y validada por los propios estudiantes.

Asimismo, el análisis muestra que la presencia de la IA introduce tensiones productivas en la actividad matemática. Las discrepancias entre las soluciones generadas por la tecnología y los modelos construidos por los estudiantes obligan a revisar supuestos, ajustar decisiones y discutir criterios de plausibilidad. Estas tensiones no interrumpen el proceso de resolución, sino que contribuyen a mantener abierta la actividad matemática y a reforzar el papel del estudiante como agente que regula el proceso.

De este modo, la complementariedad humano–IA observada en la práctica de aula se caracteriza por una colaboración funcional y dinámica, en la que la tecnología amplía las posibilidades de exploración sin sustituir las decisiones matemáticas fundamentales. Esta distribución asimétrica de funciones permite sostener una actividad matemática significativa y preservar el control epistémico del estudiante en contextos de disponibilidad de inteligencia artificial generativa.

**Prácticas de modelización y validación**

El cuarto eje de análisis se centra en las prácticas de modelización y validación que emergen durante la resolución de las tareas abiertas en contextos de disponibilidad de inteligencia artificial generativa. El análisis de la experiencia evidencia que la resolución de este tipo de tareas conduce a la construcción de modelos matemáticos basados en supuestos, los cuales son progresivamente refinados a lo largo del trabajo individual y colectivo.

Los estudiantes elaboran representaciones que articulan cantidades, relaciones y operaciones, ajustándolas en función de las decisiones adoptadas y de la coherencia con la situación planteada. La modelización no se presenta como una fase aislada, sino como un proceso iterativo en el que los modelos construidos son susceptibles de ser revisados, reformulados o descartados en función de nuevos argumentos o evidencias.

Un aspecto central de la actividad matemática observada es el papel de la validación. Los modelos no se consideran aceptables únicamente por producir un resultado numérico, sino que son sometidos a procesos de revisión en los que se evalúa su plausibilidad, consistencia interna y adecuación al contexto. Estos procesos se desarrollan tanto a partir de la discusión entre pares como del contraste con distintas fuentes de información, incluida la inteligencia artificial.

La validación adquiere así un carácter epistémico, en la medida en que implica decidir qué soluciones "tienen sentido" más allá de su corrección formal. Las discrepancias entre modelos, así como entre los resultados obtenidos por los estudiantes y los generados por la IA, actúan como catalizadores de la discusión, favoreciendo la negociación de criterios de aceptación compartidos y el examen crítico de los supuestos que sustentan las decisiones matemáticas.

En este sentido, las prácticas de modelización y validación observadas permiten cerrar el ciclo matemático de la tarea, articulando interpretación, construcción, análisis y revisión de soluciones. Este cierre refuerza el rol del estudiante como responsable último de las decisiones matemáticas y consolida una actividad matemática significativa, incluso en contextos de disponibilidad de inteligencia artificial generativa.

## Análisis o discusión de resultados

Los resultados del estudio permiten interpretar la actividad matemática observada desde una perspectiva didáctica que pone en diálogo el diseño de tareas abiertas, la regulación del uso de la inteligencia artificial y las formas de trabajo matemático que emergen en el aula. En este sentido, los hallazgos no deben entenderse como efectos directos de la presencia de la IA, sino como consecuencias de decisiones didácticas que configuran el tipo de actividad matemática que se vuelve posible en contextos post-IA.

Desde el marco propuesto por Artigue (2025), los resultados muestran un reequilibrio entre la valencia pragmática y la valencia epistémica del trabajo matemático. Si bien la inteligencia artificial facilita la obtención rápida de resultados, el diseño de las tareas abiertas y la regulación de la actividad evitan que la resolución se reduzca a una ejecución técnica. Por el contrario, los estudiantes se ven involucrados en procesos de interpretación, formulación de supuestos y validación de modelos, lo que refuerza la dimensión epistémica de la actividad matemática y preserva su sentido formativo.

En relación con la génesis instrumental, los resultados sugieren que la inteligencia artificial no se consolida como un instrumento matemático autónomo, sino como un recurso subordinado a la actividad del estudiante. La IA es utilizada para explorar o contrastar resultados, pero la información que proporciona requiere ser interpretada y validada en función de los supuestos adoptados y de la coherencia con la situación planteada. Este uso regulado se alinea con planteamientos que destacan la necesidad de una mediación didáctica para orientar la apropiación de herramientas digitales en el aprendizaje matemático (Artigue, 2025).

Desde la noción de complementariedad humano–IA (Hemmer et al., 2024), los hallazgos ponen de manifiesto una distribución funcional de roles en la que las capacidades de la tecnología y las de los estudiantes no se superponen, sino que se articulan de manera asimétrica. La inteligencia artificial amplía las posibilidades de exploración y cálculo, mientras que los estudiantes conservan el control sobre las decisiones matemáticas centrales, en particular aquellas vinculadas a la interpretación de la situación y a la validación de los resultados. Esta complementariedad no emerge de manera espontánea, sino que se construye a partir del diseño de tareas que no admiten un cierre automático.

Asimismo, las tensiones observadas entre las soluciones generadas por la IA y los modelos construidos por los estudiantes adquieren un valor didáctico relevante. Lejos de constituir un obstáculo, estas discrepancias favorecen procesos de revisión y discusión que fortalecen la actividad matemática. En este sentido, la inteligencia artificial actúa como un elemento que introduce fricción cognitiva productiva, siempre que el diseño de la tarea preserve espacios para la argumentación y la validación.

En conjunto, la discusión de los resultados permite sostener que la integración de la inteligencia artificial en matemáticas escolares no debe abordarse desde una lógica de sustitución de la actividad matemática, sino desde una reconfiguración del diseño de tareas y de las prácticas de aula. Los hallazgos refuerzan la idea de que el sentido formativo de la actividad matemática puede sostenerse en contextos post-IA cuando el diseño didáctico prioriza la toma de decisiones, la validación y la construcción de significado por parte de los estudiantes.

## Conclusiones y proyecciones

El estudio permitió analizar cómo el diseño de tareas matemáticas abiertas, articulado con una regulación didáctica explícita del uso de la inteligencia artificial, configura formas de trabajo matemático que preservan el sentido formativo de la actividad en educación secundaria. En este marco, el dispositivo COMPAS cumplió un rol central como estructura de regulación de la actividad matemática, orientando la gestión de la tarea y la interacción con la IA sin constituir una secuencia prescriptiva ni un método cerrado.

Los resultados muestran que, cuando las tareas introducen indeterminación y exigen la formulación de supuestos y la validación de resultados, los estudiantes mantienen el control epistémico del proceso matemático, incluso en contextos de disponibilidad de inteligencia artificial generativa. La implementación de COMPAS permitió regular los momentos de exploración, discusión y validación, evitando la clausura prematura de la tarea y favoreciendo una participación activa de los estudiantes a lo largo del proceso de resolución.

En particular, se evidenció que la agencia matemática del estudiante se sostiene cuando las decisiones centrales —interpretación de la situación, construcción del modelo y criterios de validación— no pueden ser delegadas a la tecnología. Asimismo, la complementariedad humano–IA observada se caracteriza por una distribución asimétrica y funcional de roles, promovida por el diseño de las tareas y por la mediación docente apoyada en COMPAS, en la que la inteligencia artificial amplía las posibilidades de exploración sin sustituir las decisiones matemáticas fundamentales.

Como limitaciones del estudio, se reconoce que el análisis se circunscribe a una experiencia de aula en un contexto específico de educación secundaria, a un conjunto acotado de tareas y a la implementación de un único dispositivo de regulación. El enfoque cualitativo adoptado no busca generalizar resultados, sino comprender en profundidad las prácticas emergentes, por lo que los hallazgos deben interpretarse a la luz de estas condiciones.

En cuanto a las proyecciones de investigación, se plantea la necesidad de explorar la implementación de COMPAS en otros niveles educativos y contenidos matemáticos, así como de analizar su articulación con modalidades de evaluación formativa en contextos con inteligencia artificial. Asimismo, futuras investigaciones podrían indagar en la incorporación de este tipo de dispositivos en la formación docente, con el fin de fortalecer el diseño de tareas y la regulación de la actividad matemática en escenarios educativos post-IA.

## Referencias